\documentclass[a4paper,12pt]{article}
\usepackage{amsmath}
\usepackage{amsfonts}
\usepackage{amssymb}
\usepackage{graphicx}
\newtheorem{theorem}{Theorem}

\newtheorem{definition}[theorem]{Definition}

\newtheorem{lemma}[theorem]{Lemma}

\newenvironment{proof}[1][Proof]{\textbf{#1.} }{\ \rule{0.5em}{0.5em}}
\begin{document}

\title{Classification of $\mathbb{Z}_{p^{k}}^{m}$ orientation preserving actions on surfaces}
\author{Antonio F.~Costa, Sergey M.~Natanzon}
\date{}

\maketitle

\emph{Dedicate to Askold Khovanskii on the occasion of his 60th
birthday.}

\vskip 2cm

\footnotetext{{\small Research partially supported by grants
RFBR-04-01-00762, INTAS 05-7805, NSh-4719.2006.1,
NWO-047.011.2004.026 (RFBR-05-02-89000-NWO-a) and Ministerio de
Ciencia y Tecnolog\'{i}a BFM 2002-04801.}}

\textbf{Abstract}. In this article we classify actions of the groups
$\mathbb{Z}_{p^{k}}^{m}$ ($p$ is a prime integer) on compact oriented surfaces.

\section{Introduction}

Finite abelian group actions on surfaces by autohomeomorphisms
constitutes a classical subject, see the articles [N], [E], [J1],
[J2], [Na1], [S], [Z], [CN]. In this article we present a direct and
complete way to deal with the topological classification of
$G$-actions, where
$G=\mathbb{Z}_{p^{k}}^{m}=\mathbb{Z}_{p^{k}}\oplus$
$\overset{m}{...}\oplus\mathbb{Z}_{p^{k}}$, $p$ is a prime integer
and $\mathbb{Z}_{p^{k}}=\mathbb{Z} /p^{k}\mathbb{Z}$.

Since $G\;$is a finite group, given an action $(\widetilde{S},f)$ of
$G$, it is possible to construct an analytic structure on $S$ such
that $f(G)$ consists of automorphisms as Riemann surface (see [K]).
Hence all the actions considered in this paper appear as
automorphism group actions of complex algebraic curves. One of the
motivations for our study is the description of the set of connected
components in the moduli space $M$ of pairs $(C,G)$, where $C$ is a
complex algebraic curve and $G\cong\mathbb{Z}_{p^{k}}^{m}$ is an
automorphisms group of $C$. According to [Na2] the description of
connected components of $M$ is reduced to the description of
topological equivalence classes of pairs $(\widetilde{S},K)$, where
$K\cong G$ is a group of autohomeomorphisms of $\widetilde{S}$. Here
we consider that $(\widetilde{S},K)$ and
$(\widetilde{S}^{\prime},K^{\prime})$ are topological equivalent if
there exist a homeomorphism $\varphi:\widetilde
{S}\rightarrow\widetilde{S}^{\prime}$ such that
$K^{\prime}=\varphi\circ K\circ\varphi^{-1}.$ Thus these topological
equivalence classes are in one-to-one correspondence with weak
equivalence classes for $G$-action on $\widetilde{S}$ in sense of
[E].

Our method for topological classification of $G$-action is a direct
generalization of the method from [CN]. It is based on topological
classification of representations $G$ in the group of all
autohomeomorphisms of $\widetilde{S}$. We prove that classes of
strong equivalence, that appear here, one-to-one correspond to
symplectic forms on $G=\mathbb{Z}_{p^{k}}^{m}$. Letter we prove,
that the classes of weak equivalence one-to-one correspond to
classes of algebraic equivalence of the symplectic forms on
$G=\mathbb{Z}_{p^{k}}^{m}$. Using this algebraic description we find
a full system of topological invariants, describing the classes of
weak equivalence. This gives very simple and natural description for
the classes of topological equivalence of the actions.

The classification of the actions of $\mathbb{Z}_{p^{k}}^{m}$ allow
to obtain strong and weak classifications for
actions by groups $\bigoplus_{i=1}^{r}\mathbb{Z}_{p_{i}^{k_i}%
}^{m_{i}}$, where $p_{i}$ are prime numbers and $p_{i}\neq p_{j}$
for $i\neq j$. Unfortunately, our methods do not work for
$p_{i}=p_{j}$. This make very involved the problem of full
classification of abelian actions.

Nevertheless, weak equivalence classes for fixed points free actions
by $G=\mathbb{Z}_{p^{2}}^m\bigoplus\mathbb{Z}_{p}^n$ was found by
Vinberg [V], using some another method. His method gives also a
simple direct formula for the number of weakly equivalence classes
for fixed points free actions by $\mathbb{Z}_{p^{k}}^{m}$.

\textbf{Acknowledgement.} Authors thank E.B.Vinberg for numerous
discussion and important consultations. This paper was written
during the authors stay at Max-Planck-Institut f\"{u}r Mathematik in
Bonn. We would like to thank this institution for its support and
hospitality.

\section{Algebraic preliminaries.}

Let us consider the standard lattice $\mathbb{Z}^{2g}=\mathbb{Z}\oplus$
$\overset{2g}{...}\oplus\mathbb{Z}$ with the standard basis $(e_{i}%
)=((0,...,1^{(i)},...0))$. We define the skew symmetric bilinear form
$(.,.):\mathbb{Z}^{2g}\times\mathbb{Z}^{2g}\rightarrow\mathbb{Z}$, by
$(e_{i},e_{j})=\delta_{i+j,2g+1}$ for $i<j$.

We consider also the group
$\mathbb{Z}_{p^{k}}^{2g}=\mathbb{Z}_{p^{k}}\oplus$
$\overset{2g}{...}\oplus\mathbb{Z}_{p^{k}}$ where $p$ is a prime and
$\mathbb{Z}_{p^{k}}=\mathbb{Z}/p^{k}\mathbb{Z}=\{\overline{0},\overline
{1},\overline{2},...,\overline{p^{k}-1}\}$. Let $\varphi:\mathbb{Z}%
^{2g}\rightarrow\mathbb{Z}_{p^{k}}^{2g}$ be the natural projection defined by
$\varphi(e_{i})=\overline{e_{i}}$, where $\overline{e_{i}}=(0,...,\overline
{1}^{(i)},...0)$. Then we have a skew symmetric bilinear form $(.,.)_{p^{k}%
}:\mathbb{Z}_{p^{k}}^{2g}\times\mathbb{Z}_{p^{k}}^{2g}\rightarrow
\mathbb{Z}_{p^{k}}$ defined by $(\overline{e_{i}},\overline{e_{j}}%
)=(\varphi(e_{i}),\varphi(e_{j}))_{p^{k}}=(e_{i},e_{j})\operatorname{mod}%
p^{k}.$ We shall denote by $\widehat{\mathbb{Z}}_{p^{k}}^{2g}$ the group
$\mathbb{Z}_{p^{k}}^{2g}$ provided with the bilinear form $(.,.)_{p^{k}}$.

Let $Sp(2g,\mathbb{Z})$ and $Sp_{p^{k}}(2g,\mathbb{Z}_{p^{k}})$ be the
subgroups of the automorphisms groups of $\mathbb{Z}^{2g}$ and $\mathbb{Z}%
_{p^{k}}^{2g}$ that preserve the bilinear forms $(.,.)$ and $(.,.)_{p^{k}}$
respectively. The natural projection $\varphi:\mathbb{Z}^{2g}\rightarrow
\mathbb{Z}_{p^{k}}^{2g}$ induces a homomorphism $\varphi_{\ast}%
:Sp(2g,\mathbb{Z})\rightarrow Sp_{p^{k}}(2g,\mathbb{Z}_{p^{k}})$ such that
$\varphi_{\ast}(f)\circ\varphi=\varphi\circ f$ for all $f\in Sp(2g,\mathbb{Z}).$

The following result is known:

\begin{lemma}
(See Lemma 4.1, pg 177 of [E]) $\varphi_{\ast}(Sp(2g,\mathbb{Z}))=Sp_{p^{k}%
}(2g,\mathbb{Z}_{p^{k}}).$
\end{lemma}

We shall also need the following result, proving in [E]

\begin{theorem}
Let $G,G^{\prime}$ be isomorphic to $\mathbb{Z}_{p^{k}}^{m}$
subgroups of $\widehat{\mathbb{Z}}_{p^{k}}^{2g}$  and
$\psi:G\rightarrow G^{\prime}$ be an isomorphism such that
$(\psi(a),\psi(b))_{p^{k}}=(a,b)_{p^{k}}$ for
all $a,b\in G$. Then there is an automorphism $\widetilde{\psi}\in Sp_{p^{k}%
}(2g,\mathbb{Z}_{p^{k}})$, such that $\psi$ is the restriction
$\widetilde{\psi}$ to $G$.
\end{theorem}

\section{Strong classification of fixed point free actions}

Let $G$ be a finite group and $\widetilde{S}$ be a closed (compact without
boundary) oriented surface. An (orientation preserving)\ action of the group
$G$ on $\widetilde{S}$ is a pair $(\widetilde{S},f)$, where $f$ is a
monomorphism of $G$ in the group of orientation preserving autohomeomorphisms
of $\widetilde{S}$. Now we consider $G=\mathbb{Z}_{p^{k}}^{m}$.

\begin{definition}
(Strong equivalence). Two actions $(\widetilde{S},f)$ and $(\widetilde
{S^{\prime}},f^{\prime})$ are called strongly equivalent if there is a
homeomorphism$,$ $\widetilde{\psi}:\widetilde{S}\rightarrow\widetilde
{S}^{\prime},$ sending the orientation of $\widetilde{S}$ to the orientation
of $\widetilde{S}^{\prime}$ and such that $f^{\prime}(h)=\widetilde{\psi}\circ
f(h)\circ\widetilde{\psi}^{-1},$ for all $h\in G.$
\end{definition}

We denote by $S=\widetilde{S}/f(G)$ and by $\varphi=\varphi(f):\widetilde
{S}\rightarrow S$ the natural projection. In this and next section we shall
consider the case when $f(h) $ has no fixed points for any $h\in G-\{id\}$, i.
e. the action of $(\widetilde{S},f)$ is fixed point free. Then the projection
$\varphi(f):\widetilde{S}\rightarrow S$ is an unbranched covering with deck
transformation group $f(G)$.

Let us consider $\pi_{1}(S)$ as the group of deck transformations of the
universal covering of $S$. Then we have:

\begin{center}
$\omega(\widetilde{S},f):\pi_{1}(S)\rightarrow\pi_{1}(S)/\pi_{1}(\widetilde
{S})=f(G)\overset{f^{-1}}{\rightarrow}G.$
\end{center}

The resulting epimorphism $\omega(\widetilde{S},f):\pi_{1}(S)\rightarrow G$ is
the monodromy epimorphism of the covering $\varphi(f):\widetilde{S}\rightarrow
S$. The epimorphism $\omega(\widetilde{S},f):\pi_{1}(S)\rightarrow G$, induces
the epimorphism $\theta(\widetilde{S},f):H_{1}(S,\mathbb{Z}_{p^{k}%
})\rightarrow G$, since $G$ is abelian.

Conversely, given an epimorphism $\theta:H_{1}(S,\mathbb{Z}_{p^{k}%
})\rightarrow G$, there is an action $(\widetilde{S},f)$ such that
$\theta=\theta(\widetilde{S},f).$ To obtain $\widetilde{S}$ it is enough to
consider the monodromy $\omega:\pi_{1}(S)\rightarrow H_{1}(S,\mathbb{Z}%
_{p^{k}})\overset{\theta}{\rightarrow}G$ and then $\widetilde{S}=U/\ker\omega
$, where $U$ is the universal covering of $S$ and the action of $G$ is given
by $G=\pi_{1}(S)/\ker\omega$.

\begin{definition}
Let $S$ and $S^{\prime}$ be two surfaces. Two epimorphisms $\theta
:H_{1}(S,\mathbb{Z}_{p^{k}})\rightarrow G$ and $\theta^{\prime}:H_{1}%
(S^{\prime},\mathbb{Z}_{p^{k}})\rightarrow G$ are called strongly equivalent
if there is an orientation preserving homeomorphism $\psi:S\rightarrow
S^{\prime}$ inducing an isomorphism $\psi_{p^{k}}:H_{1}(S,\mathbb{Z}_{p^{k}%
})\rightarrow H_{1}(S^{\prime},\mathbb{Z}_{p^{k}})$ such that $\theta
=\theta^{\prime}\circ\psi_{p^{k}}.$
\end{definition}

Using covering space theory we have that:

\begin{theorem}
Two actions $(\widetilde{S},f)$ and $(\widetilde{S}^{\prime},f^{\prime})$ are
strongly equivalent if and only if the epimorphisms $\theta(\widetilde{S},f)$
and $\theta(\widetilde{S^{\prime}},f^{\prime})$ are strongly equivalent.
\end{theorem}

\begin{definition}
Let $(\widetilde{S},f)$ be an action of $G,$ $S=\widetilde{S}/f(G)$, and
$\theta=\theta(\widetilde{S},f):H_{1}(S,\mathbb{Z}_{p^{k}})\rightarrow G$ be
the epimorphism defined by the action $(\widetilde{S},f)$. Let us consider the
spaces of homomorphisms $G^{\ast}=\{e:G\rightarrow\mathbb{Z}_{p^{k}}\}$ and
$H^{1}(S,\mathbb{Z}_{p^{k}})=\{e:H_{1}(S,\mathbb{Z}_{p^{k}})\rightarrow
\mathbb{Z}_{p^{k}}\}$. Then $\theta$ generates a monomorphism $\theta^{\ast
}=\theta^{\ast}(\widetilde{S},f):G^{\ast}\rightarrow H^{1}(S,\mathbb{Z}%
_{p^{k}}).$ The intersection form $(.,.)_{p^{k}}=(.,.)_{p^{k}}^{S}$ on
$H_{1}(S,\mathbb{Z}_{p^{2}})$ induces an isomorphism $i:H^{1}(S,\mathbb{Z}%
_{p^{k}})\rightarrow H_{1}(S,\mathbb{Z}_{p^{k}})$ defined by $(a,.)\rightarrow
a$ and a form $(.,.)_{(\widetilde{S},f)}:G^{\ast}\times G^{\ast}%
\rightarrow\mathbb{Z}_{p^{k}}$ such that $(a,b)_{(\widetilde{S},f)}%
=(i\circ\theta^{\ast}(a),i\circ\theta^{\ast}(b))_{p^{k}}.$
\end{definition}

\begin{theorem}
Two fixed point free actions $(\widetilde{S},f)$ and $(\widetilde{S}^{\prime
},f^{\prime})$ of the group $G$ are strongly equivalent if and only if
$\widetilde{S}$ and $\widetilde{S}^{\prime}$ have the same genus and
$(.,.)_{(\widetilde{S},f)}=(.,.)_{(\widetilde{S^{\prime}},f^{\prime})}.$
\end{theorem}

The proof is a direct adaptation to our case of the proof of Theorem
8 in [CN]. It is necessary only to use Lemma 1 and Theorem 2 in the
above Section instead of Theorem 1 and Theorem 3 from [CN].

\begin{theorem}
Let $(.,.):G\times G\rightarrow\mathbb{Z}_{p^{k}}$ be an skew symmetric
bilinear form. Then there exists an action $(\widetilde{S},f)$ such that
$(.,.)=(.,.)_{(\widetilde{S},f)}$ and the genus of $\widetilde{S}/f(G)$ is $g$
if and only if there exists an monomorphism $\phi:G\rightarrow\mathbb{Z}%
_{p^{k}}^{2g}$ such that $(a,b)=(\phi(a),\phi(b))_{p^{k}}$.
\end{theorem}

\begin{proof}
Let $(\widetilde{S},f)$ be an action such that $(.,.)=(.,.)_{(\widetilde
{S},f)}$. Then it is obvious by the construction of $(.,.)_{(\widetilde{S}%
,f)},$ that there exists the $\phi:G\rightarrow\mathbb{Z}_{p^{k}}^{2g}$. To
construct the action from the form and the numerical conditions, consider a
basis $(x_{i},(i=1,...,m))$ of $G$ and a surface $S$ of genus $g$. Then
$H_{1}(S,\mathbb{Z}_{p^{k}})$ has a basis $\{\chi_{i}|(i=1,...,m),\nu
_{i},(i=m+1,...,2g)\}$, such that $(\chi_{i},\chi_{j})=(x_{i},x_{j})$. We
define now the epimorphism $\theta:H_{1}(S,\mathbb{Z}_{p^{k}})\rightarrow G$,
by $\theta(\chi_{i})=x_{i},$ if $i\leq m,$ $\theta(\nu_{i})=0.$ Then the
epimorphism $\theta$ defines a homomorphism $H_{1}(S,\mathbb{Z)\rightarrow}%
H_{1}(S,\mathbb{Z}_{p^{k}})\rightarrow G$ and thus a regular covering
$\widetilde{S}\rightarrow S$ with automorphism group $G$. The action of $G$ on
$\widetilde{S}$ satisfies $(.,.)_{(\widetilde{S},f)}=(.,.)$.
\end{proof}

If $G=\bigoplus_{i=1}^{r}\mathbb{Z}_{p_{i}^{k_i}}^{m_{i}}$, where
$p_{i}$ are prime and $p_{i}\neq p_{j}$ for $i\neq j$, it is
possible to generalize our study to this more general situation. We
sketch some of the steps:

Let $(\widetilde{S},f)$ be an action of $G=\bigoplus_{i=1}^{r}\mathbb{Z}%
_{p_{i}^{k_i}}^{m_{i}}$, where $p_{i}$ are prime and $p_{i}\neq
p_{j}$ for $i\neq j$. Put $n=p_{1}^{k_{1}}...p_{r}^{k_{r}}$. Put
$S=\widetilde{S}/f(G)$, and
$\theta=\theta(\widetilde{S},f):H_{1}(S,\mathbb{Z}_{n})\rightarrow
G$ be the epimorphism defined by the action $(\widetilde{S},f)$. Let
us consider the spaces of homomorphisms
$G^{\ast}=\{e:G\rightarrow\mathbb{Z}_{n}\}$ and
$H^{1}(S,\mathbb{Z}_{n})=\{e:H_{1}(S,\mathbb{Z}_{n})\rightarrow\mathbb{Z}%
_{n}\}$. Then $\theta$ generates a monomorphism $\theta^{\ast}=\theta^{\ast
}(\widetilde{S},f):G^{\ast}\rightarrow H^{1}(S,\mathbb{Z}_{n}).$ For each
$i\in\{1,...,r\}$, the intersection form $(.,.)_{n}^{S}$ on $H_{1}%
(S,\mathbb{Z}_{n})$ induces a form $(.,.)_{(\widetilde{S},f)}^{(i)}%
:\mathbb{Z}_{p_{i}^{k_i}}^{m_{i}}\times\mathbb{Z}_{p_{i}^{k_i}}^{m_{i}}%
\rightarrow\mathbb{Z}_{p_{i}^{k_i}}$. The set of forms
$\{(.,.)_{(\widetilde {S},f)}^{(i)}\}$ and the genus of the surface
$\widetilde{S}$ gives now the complete set of invariants for the
strong classification.

\section{Weak classification of fixed point free actions}

\begin{definition}
(Weak equivalence) Let $(\widetilde{S},f)$ and $(\widetilde{S^{\prime}%
},f^{\prime})$ be two actions of a group $G$. We shall say that $(\widetilde
{S},f)$ and $(\widetilde{S^{\prime}},f^{\prime})$ are weakly equivalent if
there is a homeomorphism $\widetilde{\psi}:\widetilde{S}\rightarrow
\widetilde{S}^{\prime}$ and an automorphism $\alpha\in Aut(G)$ such that
$f^{\prime}\circ\alpha(h)=\widetilde{\psi}\circ f(h)\circ\widetilde{\psi}%
^{-1}$, $h\in G.$
\end{definition}

The map $v\mapsto p^{k-1}v$ generate the natural epimorphism to the
vector space $\mu: \mathbb{Z}_{p^k}^m\rightarrow \mathbb{Z}_p^m$.
For an action $(\widetilde{S},f)$ let us put
$q(\widetilde{S},f)=(q_1(\widetilde{S},f),...,q_k(\widetilde{S},f))$,
where $q_i(\widetilde{S},f)$ is the dimension of the vector space
$\mu(\{h\in G^{\ast}|(h,G^{\ast })_{(\widetilde{S},f)}\equiv 0(\mod
p^i)\})\subset\mu(G^{\ast })\cong\mathbb{Z}_p^m$. It is follow from
[B, Ch IX,5,Theorem 1], that  $0\leq q_k\leq...\leq q_1\leq m$ and
$q_i\equiv m\mod2$.

{\ }

The next Theorem solves the problem of weak classification of actions of
$\mathbb{Z}_{p^{k}}^{m}$on surfaces:

\begin{theorem}
Let $(\widetilde{S},f)$ and $(\widetilde{S^{\prime}},f^{\prime})$ be
two fixed point free actions of a group $G=\mathbb{Z}_{p^{k}}^{m}.$
Then the actions $(\widetilde {S},f)$ and $(\widetilde{S^{\prime}},
f^{\prime})$ are weakly equivalent if and only if $\widetilde{S}$
and $\widetilde{S}^{\prime}$ have the same genus and
$q(\widetilde{S},f)=q(\widetilde{S^{\prime}},f^{\prime}).$
\end{theorem}

\begin{proof}
Let $S=\widetilde{S}/f(G)$ and $S^{\prime}=\widetilde{S^{\prime}}/f^{\prime
}(G)$ have the same genus $g$ and $q(\widetilde{S},f)=q(\widetilde{S^{\prime}%
},f^{\prime}).$ Let $\theta^{\ast}(\widetilde{S},f)$ and $\theta^{\ast
}(\widetilde{S}^{\prime},f^{\prime})$ be the epimorphisms defined by the two
actions. Consider the image $\widetilde{G}$ of $G^{\ast}$ in $H_{1}%
(S,\mathbb{Z}_{p^{k}})$ by $\theta^{\ast}(\widetilde{S},f)$ and the image
$\widetilde{G}^{\prime}$ of $G^{\ast}$ in $H_{1}(S^{\prime},\mathbb{Z}_{p^{k}%
})$ by $\theta^{\ast}(\widetilde{S}^{\prime},f^{\prime})$. Then
there exists an isomorphism
$\psi:\widetilde{G^{\prime}}\rightarrow\widetilde{G}$ such that
$(\psi(a),\psi(b))_{(\widetilde{S^{\prime}},f^{\prime})}=(a,b)_{(\widetilde
{S},f)}$. It follows from Theorem 2, Lemma 1 and [B, Ch IX,5,Theorem
1 ], that there exists an isomorphism
$\tilde{\psi}:H^{1}(S^{\prime},\mathbb{Z})\rightarrow
H^{1}(S,\mathbb{Z})$ giving by restriction $\psi$ and sending the
intersection form of $H^{1}(S^{\prime},\mathbb{Z})$ to the
intersection form of $H^{1}(S,\mathbb{Z})$.

Let $(.,.)_{(\widetilde{S},f)}$ and $(.,.)_{(\widetilde{S^{\prime}},f^{\prime
})}$ be the skew symmetric bilinear forms induced by the two actions. Assume
that $(.,.)_{(\widetilde{S},f)}$ and $(.,.)_{(\widetilde{S^{\prime}}%
,f^{\prime})}$ have the same signature then there exists an
isomorphism $\psi:\widetilde{G}^{\prime}\rightarrow\widetilde{G}$
such that $(\psi
(a),\psi(b))_{(\widetilde{S}^{\prime},f^{\prime})}=(a,b)_{(\widetilde{S},f)}.$
By [MKS, page 178], there exists a homeomorphism
$\varphi:S\rightarrow S^{\prime}$ inducing $\widetilde{\psi}$ on
cohomology. Then by theorem 7, the actions $(\widetilde{S},f)$ and
$(\widetilde{S},\varphi^{-1}\circ f^{\prime }\circ\varphi)$ are
strongly equivalent. The isomorphism $\psi$, defines an automorphism
of $G$ giving the weak equivalence between $(\widetilde{S},f)$ and
$(\widetilde{S^{\prime}},f^{\prime}).$
\end{proof}

Using a modification of the method from [CN, Theorem 9]  we can
construct a fixed point free actions $(\widetilde{S},f)$ with
$q(\widetilde{S},f)=q$ for any set of integer numbers
$q=(q_1,...,q_k)$, such that $0\leq q_k\leq...\leq q_1\leq m$ and
$q_i\equiv m\mod2$. Moreover, an action $(\widetilde{S},f)$ with
these invariants exists if and only if the genus of the surface
$\widetilde{S}$ is not less that $p^{km}(\frac{1}{2}(m+q_1)-1)+1$.
This gives a full description for the set of weakly equivalence
classes.

As in the above Section it is possible to extend the results of this
section to the case
$G=\bigoplus_{j=1}^{r}\mathbb{Z}_{p_{j}^{k_j}}^{m_{j}}$, where
$p_{i}$ are prime and $p_{i}\neq p_{j}$ for $i\neq j$. In this case
the invariant $q(\widetilde{S},f)$ must be substituted to the set of
invariants $q_i^j(\widetilde{S},f)$, corresponding to the
$\mathbb{Z}_{p_{j}^{k_j}}^{m_{j}}$ actions.

\section{Classification of actions with elements having fixed points.}

Consider now an arbitrary effective action $(\tilde{S},f)$ of $G=\mathbb{Z}%
_{p^{k}}^{m}$. Let $B\subset S=\tilde{S}/f(G)$ be the set of all critical
values of the natural projection $\varphi: \tilde{S}\rightarrow S$. The
construction from section 3 define now the epimorphism $\theta(\widetilde
{S},f):H_{1}(S-B,\mathbb{Z}_{p^{k}})\rightarrow G$.

The images by $\theta$ of elements from $H_{1}(S-B,\mathbb{Z}_{p^{k}})$,
surrounding points $b\in B$ generate the subgroup $G_{fix}\subset G$. This
group contains all $h\in G$, such that $f(h)$ has fixed points, and it is
generated by such $h^{\prime}s$. Consider the function $l=l_{(\widetilde
{S},f)}:G\rightarrow\mathbb{Z}_{\geq0}=\{n\in\mathbb{Z}|n\geq0\}$, where
$l(h)$ is the number of elements $b\in B$, such that $\theta(b)=h$. The
functions with this properties will be called \textit{characteristic functions
of $(\widetilde{S},f)$}. It is obvious that $\sum_{h\in G}l(h)h=0$. Two
functions $l,l^{\prime}:G\rightarrow\mathbb{Z}_{\geq0}$ we call
\textit{equivalent} if there exists a automorphism $\alpha\in Aut(G)$ such
that $l^{\prime}(h)=l(\alpha(h))$.

Assume now $G_{free}=G/G_{fix}\cong\mathbb{Z}_{p^{k}}^{m}$ and we consider the
epimorphism $\vartheta:H_{1}(S,\mathbb{Z}_{p^{k}})\rightarrow G_{free}$,
defined by $H_{1}(S,\mathbb{Z}_{p^{k}})\rightarrow H_{1}(S-B,\mathbb{Z}%
_{p^{k}})/G_{fix}\rightarrow G/G_{fix}=G_{free}$. In fact, the
epimorphism $\vartheta$ defined by fixed point-free action of
$G_{free}$ on $\tilde {S}/G_{fix}$. Thus the construction from
section 3 and 4 gives the skew symmetric bilinear form
$(.,.)_{(\widetilde{S},f)}:G_{free}^{\ast}\times G_{free}^{\ast
}\rightarrow\mathbb{Z}_{p^{k}}$ and the set of numbers
$q(\widetilde{S},f)$ for this form.

Using the same arguments that in Theorem 13,14,15 of [CN] we prove.

\begin{theorem}
Two actions $(\widetilde{S},f)$ and $(\widetilde{S}^{\prime},f^{\prime})$ of
the group $G=\mathbb{Z}_{p^{k}}^{m}$ with $G_{free}\cong\mathbb{Z}_{p^{k}}%
^{n}$ are strongly equivalent if and only if $\widetilde{S}$ and
$\widetilde{S}^{\prime}$ have the same genus, $l_{(\widetilde{S}%
,f)}=l_{(\widetilde{S^{\prime}},f^{\prime})}$ and $(.,.)_{(\widetilde{S}%
,f)}=(.,.)_{(\widetilde{S^{\prime}},f^{\prime})}.$
\end{theorem}

\begin{theorem}
Two actions $(\widetilde{S},f)$ and
$(\widetilde{S}^{\prime},f^{\prime})$ of
the group $G=\mathbb{Z}_{p^{k}}^{m}$ with $G_{free}\cong\mathbb{Z}_{p^{k}}%
^{n}$ are strongly equivalent if and only if $\widetilde{S}$ and
$\widetilde{S}^{\prime}$ have the same genus, $l_{(\widetilde{S}%
,f)}=l_{(\widetilde{S^{\prime}},f^{\prime})}$ and
$q(\widetilde{S},f)=q(\widetilde{S}',f').$
\end{theorem}

\begin{center}
{\LARGE References}
\end{center}

[B] N. Bourbaki, \'{E}l\'{e}ments de math\'{e}matique. Premi\`{e}re partie:
Les structures fondamentales de l'analyse. Livre II: Alg\`{e}bre. Chapitre 9:
Formes sesquilin\'{e}aires et formes quadratiques. (French) Actualit\'{e}s
Sci. Ind. no. 1272 Hermann, Paris 1959 211 pp.

[CN] A. F. Costa and S. M. Natanzon, Topological Classification of
$\mathbb{Z}_{p}^{m}$ Actions on Surfaces, \textit{Michigan Math. J}\emph{.}
\textbf{50} (2002)\ 451-460.

[E] A. L. Edmonds, Surface symmetry I, \textit{Michigan Math. J.},
\textbf{29}, (1982), 171-183.

[J1] S. A. Jassim, Finite abelian coverings, \textit{Glasgow Math. J}.,
\textbf{25}, (1984), 207-218.

[J2] S. A. Jassim, Finite abelian actions on surfaces, \textit{Glasgow Math.
J}., \textbf{35}, (1993), 225-234.

[MKS] W. Magnus, A. Karras and D. Solitar, Combinatorial group theory, Dover,
New York, 1976.

[K] B.Kerekjarto, Vorlesungen uber Topologie. I. Flachentopologie. Berlin:
Springer-Verlag, 1923.

[Na1] S. M. Natanzon, Moduli spaces of complex algebraic curves with
isomorphic to $(\mathbb{Z}/2\mathbb{Z})^{m}$ groups of automorphisms,
\textit{Differential Geometry and its Applications} \textbf{5}, (1995) 1-11.

[Na2] S. M. Natanzon, Moduli of Riemann surfaces, Hurwitz-type spaces, and
their superanalogues, \textit{Russian Math. Surveys}, \textbf{54}:1 (1999) 61-117.

[N] J. Nielsen, Die Struktur periodischer Transformationen von Fl\"{a}chen,
\textit{Danske Vid. Selsk. Mat.-Fys. Medd}. \textbf{15} (1937), 1-77.

[S] P. A. Smith, Abelian actions of periodic maps on compact surfaces,
\textit{Michigan Math. J.} \textbf{14} (1967), 257-275.

[V] E. B. Vinberg, On abelian covering of surfaces
(\textit{presented to Michigan Math. J}).

[Z] B. Zimmermann, Finite abelian group actions on surfaces, \textit{Yocohama
Mathematical Journal,} \textbf{38}, (1990) 13-21.

\bigskip

\bigskip

{\small Antonio F. Costa}

{\small Departamento de Matem\'{a}ticas Fundamentales,UNED}

28040-{\small Madrid, Spain.}

{\small e-mail:\ acosta@mat.uned.es}

\bigskip

{\small Sergey M. Natanzon}

{\small Moscow State University, Independent University of Moscow, ITEP}

{\small Moscow, Russia.}

{\small e-mail: natanzon@mccme.ru}

\
\end{document}